\documentclass [12pt,a4paper]{amsart}
\usepackage{cases}
\usepackage{amscd}
\usepackage{epsfig}
\usepackage{graphicx}
\usepackage{verbatim}
\usepackage{amssymb}
\usepackage{amsfonts}
\usepackage{amsbsy}
\usepackage{latexsym,amsmath,amsthm,amsfonts,mathrsfs}
\usepackage{amsthm}
\usepackage{amstext}
\usepackage{amsopn}
\usepackage[dvips]{color}
\usepackage{indentfirst}

\newtheorem{thm}{Theorem}[section]
\newtheorem{defi}[thm]{Definition}

\newtheorem{lem}[thm]{Lemma}

\newtheorem{rem}[thm]{Remark}
\newtheorem{conj}[thm]{Conjecture}

\author[X. Huang]{Xueping Huang}
\address{Department of Mathematics\\ University of Bielefeld, 33501
Bielefeld , Germany}
\thanks{Research supported by Project CRC701}
\email{xhuang1@math.uni-bielefeld.de}

\title[]{Stochastic incompleteness for graphs and weak Omori-Yau maximum principle}
\begin{document}

\begin{abstract}
We prove an analogue of the weak Omori-Yau maximum principle and
Khas'minskii's criterion for graphs in the general setting of Keller
and Lenz. Our approach naturally gives the stability of stochastic
incompleteness under certain surgeries of graphs. It allows to
develop a unified approach to all known criteria of stochastic
completeness/incompleteness, as well as to obtain new criteria.
\end{abstract}

\date{\today}


\maketitle
\section*{Introduction}
Recently, 
Wojciechowski \cite{WojThesis}, \cite{WojSurvey}, \cite{WojIndiana}
and Weber \cite{Weber} independently studied the problem of
stochastic completeness for the following Laplace operator (Weber
calls it ``physical Laplacian'') on a locally finite, connected,
undirected graph $(V, E)$ :
\begin{equation*}
  \Delta f(x)=\sum_{y\in V, y\sim x}(f(x)-f(y))
\end{equation*}
where $V$ is the set of vertices and $E$ is the set of edges, and
$y\sim x$ means that $(x,y)\in E$. See also \cite{Jorg-Pearse09} for
some remarks. 
Essential
self-adjointness of $\Delta$ has been shown by several authors
independently, see \cite{Jorgensen08, Masamune, Weber, WojThesis}.
The corresponding heat semigroup can be constructed as $P_t
=\exp(-t\Delta)$. This semigroup determines a continuous time random
walk on $V$, that is stochastically complete provided $P_t 1 = 1$,
and incomplete otherwise. The latter can occur due to a very fast
escape rate so that the random walk reaches infinity in finite time.
This phenomenon in the setting of Brownian motions on manifolds was
first observed by Azencott \cite{Azencott} (see also the survey
\cite{Grigor}).

The study of continuous time Markov chains has a long history, see
for example the work of Feller \cite{Feller1} \cite{Feller2} and
Reuter \cite{Reuter}. However, in the analytic study of random walks
on graphs, the phenomenon of stochastic incompleteness has remained
unnoticed until recently, perhaps because most attention was given
to the normalized (combinatorial) Laplace operator
$$\tilde{\Delta}f(x)=\frac{1}{\deg x}\sum_{y\in V, y\sim x}(f(x)-f(y)),$$
where $\deg x$ is the number of neighbors of $x$ in the graph. The
corresponding heat semigroup of $\tilde{\Delta}$ is always
stochastically complete, which is a consequence of the boundedness
of $\tilde{\Delta}$ in $L^2$. Following the classical approach
\cite{Khas} to the stochastic completeness in the framework of
continuous spaces, Wojciechowski showed the equivalence of
stochastic incompleteness and the existence of the so-called
$\lambda$-(sub)harmonic functions. Wojciechowski used this
equivalence to obtain many interesting sufficient conditions for
stochastic completeness and incompleteness.

Weber \cite{Weber} followed another approach via bounded solutions
of the heat equation and discovered an interesting curvature-type
criterion.

Keller and Lenz \cite{Keller-Lenz} have extended their work to a
general setting, namely regular Dirichlet forms on discrete
countable sets. In this setting the graphs are not necessarily
locally finite and have general weight functions both for vertices
and edges.

In this note we adopt an alternative approach to stochastic
completeness on graphs. A cornerstone of this approach is Theorem
\ref{equiv-thm}, where we prove that the stochastic completeness is
equivalent to a discrete analogue of the weak Omori-Yau maximum
principle. The latter notion was introduced by Pigola, Rigoli, and
Setti \cite{PRSProc}, \cite{PRSSurvey}, where they proved the
aforementioned equivalence in the setting of manifolds and gave many
applications. For the original form of Omori-Yau maximum principle,
see \cite{Omori}, \cite{Yau}.

We use the weak Omori-Yau maximum principle and its consequence, a
discrete Khas'minskii's criterion, to develop a unified approach to
all known criteria of stochastic completeness/incompleteness, as
well as to obtain new criteria. For example, in Theorem
\ref{curvature-thm}, we establish an improvement of the
curvature-type criterion in \cite{Weber}. Together with Lemma
\ref{comparison-lem}, the weak Omori-Yau maximum principle also
easily gives stability results for stochastic incompleteness.  For
example, the subgraph of a stochastically incomplete graph, which
consists of vertices with weighted degrees larger than some
constant, is stochastically incomplete as well (Theorem
\ref{stability-sub}). Due to special features of the graph case,
some results are new and some are stronger than their manifold
relatives. Part (3) of Lemma \ref{comparison-lem} and Theorem
\ref{global-degree} have no analogues for manifolds to the author's
knowledge. Our version of Khas'minskii's criterion Theorem
\ref{Khas-thm} is stronger than a direct generalization of the
manifold case (Theorem \ref{weak-Khas}).

The paper is organized as follows. In Section 1, we introduce the
framework of Keller and Lenz as our starting point. The weak
Omori-Yau maximum principle for graphs is proved in Section 2
together with a useful Lemma \ref{comparison-lem}. Then
Khas'minskii's criterion is established in Section 3. Section 4 is
devoted to the stability of stochastic incompleteness under certain
surgeries of graphs. In Section 5, we concentrate on the special
case of physical Laplacian and show how the weak maximum principle
and Khas'minskii's criterion are applied. In the last section, we
present some open questions and further developments.

\section{Foundations}
We generally follow the framework set up in \cite{Keller-Lenz}
except that we don't include killing terms here. Consider a triple
$(V, b, \mu)$ where $V$ is a discrete countably infinite set, $\mu$
is a measure on $V$ with full support, and $b: V\times V\rightarrow
[0,+\infty)$ satisfies:

 \begin{enumerate}
 \item $b(x,x)=0$;
 \item $b(x,y)=b(y,x)$;
 \item $\sum_{y\in V}b(x,y)<+\infty$.
\end{enumerate}

The triple $(V, b, \mu)$ will be called a (weighted) graph, and
sometimes we abuse the notation and denote a graph simply by $V$.
 We call the quantity $$\text{\rm{Deg}}(x):= \frac{1}{\mu(x)}\sum_{y\in
 V}b(x,y)$$
the weighted degree of $x\in V$. For example, for the physical
Laplacian, $\mu(x)\equiv1$ and the weighted degree
$\text{\rm{Deg}}(x)$ coincides with the usual degree $\deg x$.
However, for the combinatorial Laplacian, $\mu(x)=\deg(x)$ and hence
$\text{\rm{Deg}}(x)\equiv 1$.

The couple $(V,\mu)$ forms a measure space. Then the real function
spaces $L^p(V,\mu), 0<p<\infty$ are naturally defined as
$$\{u:V\rightarrow \mathbb{R}: \sum_{x\in V}\mu(x)|u(x)|^p
<\infty\}$$ and $L^{\infty}(V,\mu)$ is simply the space of bounded
functions on $V$.

A \emph{formal} Laplacian $\Delta$:
  $$\Delta u(x)=\frac{1}{\mu(x)}\sum_{y}b(x,y)(u(x)-u(y))$$
is introduced on the domain
  $$F=\{u:V\rightarrow \mathbb{R}:\forall x\in V, \sum_{y}b(x,y)|u(y)|<\infty\}.$$
An obvious fact is that $L^{\infty}(V,\mu)\subseteq F$.

A quadratic form $Q$ can be defined on the space of finitely
supported functions $C_c (V)$ as
\begin{equation*}
 Q(u)=\frac{1}{2}\sum_{x,y\in V}b(x,y)(u(x)-u(y))^2.
\end{equation*}
It is closable and its closure is a regular Dirichlet form which we
also denote by $Q$. This is a nonlocal Dirichlet form in general.
The semigroup $P_t$ corresponding to the Dirichlet form $Q$ on
$L^2(V,\mu)$ can be extended to all $L^p (V,\mu), p\in[1,\infty]$,
and the associated generators are certain restrictions of the
\emph{formal} Laplacian $\Delta$. We abuse the notation and denote
all these operators by $\Delta$. The explicit domains of these
generators are irrelevant to the problem of stochastic completeness
(for details see \cite{Keller-Lenz}). For the general theory of
Dirichlet forms and semigroups, we
refer to \cite{Fukushima-O-T, MaRockner, Davies}. 

The following theorem about stochastic incompleteness is classical
for the manifold case \cite{Grigor} and is proven independently by
Wojciechowski and Weber for the graph case (physical Laplacian).
Keller and Lenz \cite{Keller-Lenz, Keller-Lenz-10} prove it in the
general setting:
\begin{thm}
  The following statements are equivalent:
\begin{enumerate}
\item For some $t_0>0$, some $x_0 \in V$, $P_{t_0}1(x_0)<1$.

 \item For every $\lambda>0$, there exists a nonconstant, nonnegative, bounded function $v$ on $V$ such that $\Delta v +\lambda v
  =0$. Such a function $v$ is called a $\lambda$-harmonic
  function.

  \item For every (or, equivalently, for some) $\lambda>0$, there exists a nonconstant, nonnegative, bounded function $v$ on $V$ such that $\Delta v +\lambda v
  \leq 0$. $v$ is called a $\lambda$-subharmonic
  function.

  \item There exists a nonconstant, nonnegative, bounded solution
  to the Cauchy problem
\begin{equation*}
\left\{ \begin{aligned}
         &\Delta u(x,t)+\frac{\partial u}{\partial t}(x,t)=0, \text{for all ~~} x\in V,\text{all~~} t\geq 0 \\
          &u(\cdot,0)=0
                          \end{aligned} \right.
                          \end{equation*}
\end{enumerate}
\end{thm}

\begin{rem}\rm{
In the context of Markov chains, Feller \cite{Feller1, Feller2} and
Reuter \cite{Reuter} also cover parts of this result from the
probabilistic point of view.}
\end{rem}

A graph is said to be stochastically incomplete if any one of these
four conditions holds. Otherwise it is stochastically complete.

\section{Weak Omori-Yau maximum principle}
From now on, we will denote the supremum of a function $f$ by $f^*$.
\begin{defi}
  A graph $(V, b, \mu)$ is said to satisfy the weak Omori-Yau maximum principle if
  for every nonnegative function $f$ on $V$ with $f^* =\sup_{V}f<+\infty$ and for
  every $\alpha>0$,
  $$\sup_{\Omega_{\alpha}}\Delta f \geq 0,$$
  where $$\Omega_{\alpha}=\{x\in V: f(x)>f^* -\alpha\}.$$
\end{defi}

It was first noticed by Pigola, Rigoli and Setti \cite{PRSProc} that
in fact a smooth, connected, non-compact, Riemannian manifold
satisfies the weak Omori-Yau maximum principle if and only if the
semigroup generated by the Laplace-Beltrami operator on it is
stochastic complete. It is somewhat surprising that this also holds
in the graph case although we are dealing with nonlocal operators
here.
\begin{thm}\label{equiv-thm} A graph satisfies the weak Omori-Yau maximum principle
if and only if it is stochastically complete.
\end{thm}
\begin{proof}
Assume that the weak maximum principle holds but the graph is
stochastically incomplete. Then there exists a bounded,
non-negative, nonconstant solution $f$ of the equation $\Delta f
+\lambda f=0$ for some $\lambda>0$. Choosing $\alpha =\frac{f^*} {2}
>0$, we have
$$\sup_{\Omega_{\alpha}}\Delta f=\sup_{\Omega_{\alpha}}-\lambda f\leq -\lambda \frac{f^*}{2}<0$$
which is a contradiction.

Conversely, if $V$ is stochastically complete but the weak maximum
principle does not hold, there exists a nonnegative function $f$ on
$V$ with $f^* <+\infty$ and some $\alpha>0$ and $c>0$ such that
$$\sup_{\Omega_{\alpha}}\Delta f<-2c.$$

Define $$f_{\alpha}=(f+\alpha -f^* )_{+},$$ which is obviously
nonconstant, nonnegative and bounded. Setting
$\lambda=\frac{c}{\alpha}$, we claim that
$$\Delta f_{\alpha} +\lambda f_{\alpha}
  \leq 0,$$
which implies stochastic incompleteness and leads to a
contradiction.

For $x\in\Omega_{\alpha}^{c}$, $f_{\alpha}(x)=0$, so the claim is
trivially true.

For $x\in\Omega_{\alpha}$, we have $$\lambda f_{\alpha}(x)\leq
\lambda\alpha=c,$$ and
$$f_{\alpha}(x)-f_{\alpha}(y)=f(x)-f^{*}+\alpha-f_{\alpha}(y)\leq f(x)-f(y).$$

Hence
\begin{equation*}
\begin{split}
\Delta f_{\alpha}(x) +\lambda
f_{\alpha}(x)&=\frac{1}{\mu(x)}\sum_{y}b(x,y)(f_{\alpha}(x)-f_{\alpha}(y))+\lambda
f_{\alpha}(x)
 \\
&\leq \frac{1}{\mu(x)}\sum_{y}b(x,y)(f(x)-f(y))
+c\\
&=\Delta f(x) +c\leq -c.
 \end{split}
 \end{equation*}
\end{proof}

Now we know that a graph is stochastically incomplete if and only if
there exist a nonnegative function $f$ on $V$ with $f^* <+\infty$
and some $\alpha>0$ and $c>0$ such that
$$\sup_{\Omega_{\alpha}}\Delta f<-c.$$

The following lemma describes some elementary properties of a
function $f$ that violates the weak maximum principle.
\begin{lem}\label{comparison-lem}
Suppose that $(V,b,\mu)$ is stochastically incomplete. Let $f$ be a
nonnegative function on $V$ such that $f^* <+\infty$ and for some
$\alpha>0$ and $c>0$,
$$\sup_{\Omega_{\alpha}}\Delta f<-c.$$

Let $\alpha'=\min\{\alpha, c\}$. Then the following is true.
\begin{enumerate}
 \item $f$ cannot attain its supremum $f^*$ on $V$, and in particular, is nonconstant;
 \item $\sup_{\Omega_{\alpha'}}\Delta f<-\alpha'$;
 \item for every $n\geq 1$, and every $x\in\Omega_{\frac{\alpha'}{n}}$,
$$\text{\rm{Deg}}(x)=\frac{1}{\mu(x)}\sum_{y}b(x,y)>n.$$ In other words,
$$\Omega_{\frac{\alpha'}{n}}\subseteq \{x\in V: \text{\rm{Deg}}(x)> n\}.$$
\end{enumerate}
\end{lem}
\begin{proof}
  (1) Suppose that there exists $x_0 \in V$ such that $f(x_0 )=f^*$. In particular, $x_0 \in \Omega_{\alpha}$. We have
that
$$\frac{1}{\mu(x_0)}\sum_{y\in V
}b(x_0,y)(f(y)-f(x_0))=-\Delta f(x_0 )>c>0.$$ Thus there exists
$y\in V$ such that $f(y)>f(x_0)$, a contradiction.

(2) Since $\alpha'\leq \alpha$, we have
$\Omega_{\alpha'}\subseteq\Omega_{\alpha}$. So
$$\sup_{\Omega_{\alpha'}}\Delta f\leq\sup_{\Omega_{\alpha}}\Delta
f<-c\leq-\alpha'.$$

(3) For $x\in\Omega_{\frac{\alpha'}{n}}$, set
$$l=\frac{1}{\mu(x)}\sum_{y:f(y)>f(x)}b(x,y),$$ we have
$$\alpha'<-\Delta
f(x)\leq\frac{1}{\mu(x)}\sum_{y: f(y)>f(x)}b(x,y)(f(y)-f(x))\leq
\frac{l\alpha'}{n}.$$

Therefore $l>n$ and, in particular, $\text{\rm{Deg}}(x)>n$ for all
$x\in \Omega_{\frac{\alpha'}{n}}$.
\end{proof}

\begin{rem}\rm{
Part (3) of this lemma gives a control of ``directions of increase"
of the function violating the weak maximum principle. An immediate
consequence is that a stochastically incomplete graph necessarily
has unbounded weighted degree. In particular, the semigroup
corresponding to the combinatorial Laplacian on a graph is
stochastically complete. This is a result of \cite{Dodziuk1, Dodziuk2}.}
\end{rem}

Stochastic incompleteness is a global property while the weighted
degree function is a local quantity. We can define a ``global
weighted degree function'' in an iterative way.
\begin{lem}\label{nonincreasing}
Fix a non-decreasing sequence $\Theta=\{a_k\}_{k\geq 0}$ of
nonnegative real numbers. We use the convention that
$$\sum_{y, y\in\emptyset}b(x,y)=0.$$
 For $x\in V$ and $k\in \mathbb{N}$, define
$$\text{\rm{Deg}}_{\Theta, 0}(x)=\text{\rm{Deg}}(x),$$ and
$$\text{\rm{Deg}}_{\Theta, k+1}(x)=\frac{1}{\mu(x)}\sum_{y,\text{\rm{Deg}}_{\Theta, k}(y)>a_k}b(x,y).$$
Then for any $x\in V$, $\{\text{\rm{Deg}}_{\Theta, k}(x)\}_{k\geq
0}$ forms a non-increasing, nonnegative sequence. In particular,
$$\text{\rm{Deg}}_{\Theta, \infty}(x)=\lim_{k\rightarrow\infty}\text{\rm{Deg}}_{\Theta, k}(x)$$
exists for all $x\in V$.
\end{lem}
\begin{proof}
The sequence $\{\text{\rm{Deg}}_{\Theta, k}(x)\}_{k\geq 0}$
obviously has nonnegative entries. We only need to prove that for
any $k\geq 0$,
$$\text{\rm{Deg}}_{\Theta, k+1}(x)\leq\text{\rm{Deg}}_{\Theta, k}(x).$$
For $k=0$, we have
$$\text{\rm{Deg}}_{\Theta, 1}(x)=\frac{1}{\mu(x)}\sum_{y,\text{\rm{Deg}}(y)>a_0}b(x,y)\leq\frac{1}{\mu(x)}\sum_{y}b(x,y)
=\text{\rm{Deg}}_{\Theta, 0}(x).$$

Assume that the assertion holds for $k=n-1\geq0$, that is
$$\text{\rm{Deg}}_{\Theta, n}(x)\leq\text{\rm{Deg}}_{\Theta, n-1}(x).$$
Since $a_n \geq a_{n-1}$, we see that for $k=n$,
\begin{equation*}
  \begin{split}
\text{\rm{Deg}}_{\Theta,
n}(x)&=\frac{1}{\mu(x)}\sum_{y,\text{\rm{Deg}}_{\Theta,
n-1}(y)>a_{n-1}}b(x,y)\\
&\geq\frac{1}{\mu(x)}\sum_{y,\text{\rm{Deg}}_{\Theta,
n-1}(y)>a_{n}}b(x,y)\\&\geq\frac{1}{\mu(x)}\sum_{y,\text{\rm{Deg}}_{\Theta,
n}(y)>a_{n}}b(x,y)\\&=\text{\rm{Deg}}_{\Theta, n+1}(x).
  \end{split}
\end{equation*}
The assertion follows by induction.
\end{proof}
\begin{defi}
We call $\text{\rm{Deg}}_{\Theta, \infty}(x)$ the global weighted
degree of $x$ with respect to the sequence
 $\Theta$. For the special case when $a_k \equiv n$, we denote $\text{\rm{Deg}}_{\Theta,
\infty}(x)$ by $\text{\rm{Deg}}_{n, \infty}(x)$ and call it the
global weighted degree of $x$ with parameter $n$.
\end{defi}
\begin{lem}\label{degm<degn}
For $m>n\geq 1$, $k\in \mathbb{N}$, the following holds for any
$x\in V$,
$$\text{\rm{Deg}}_{n,k}(x)\geq \text{\rm{Deg}}_{m,k}(x).$$
In particular, for any $x\in V$,
$$\text{\rm{Deg}}_{n,\infty}(x)\geq
\text{\rm{Deg}}_{m,\infty}(x).$$
\end{lem}
\begin{proof}
This can be proven by an induction procedure similar to the proof of
Lemma \ref{nonincreasing}. The $k=0$ case is obvious as
$$\text{\rm{Deg}}_{n,0}(x)=\text{\rm{Deg}}(x)=\text{\rm{Deg}}_{m,0}(x).$$
Assume that
$$\text{\rm{Deg}}_{n,k}(x)\geq \text{\rm{Deg}}_{m,k}(x).$$
Then we have
\begin{equation*}
  \begin{split}
\text{\rm{Deg}}_{n,k+1}(x)&=\frac{1}{\mu(x)}\sum_{y,\text{\rm{Deg}}_{n,
k}(y)>n}b(x,y)\\&\geq\frac{1}{\mu(x)}\sum_{y,\text{\rm{Deg}}_{m,
k}(y)>n}b(x,y)\\&\geq\frac{1}{\mu(x)}\sum_{y,\text{\rm{Deg}}_{m,
k}(y)>m}b(x,y)\\&=\text{\rm{Deg}}_{m,k+1}(x).
  \end{split}
\end{equation*}
\end{proof}

The notion of the global weighted degree function allows us to
improve Lemma \ref{comparison-lem} as follows.
\begin{thm}\label{global-degree}
Suppose that $(V,b,\mu)$ is stochastically incomplete. Let $f$ be a
nonnegative function on $V$ such that $f^* <+\infty$ and for some
$\alpha>0$,
$$\sup_{\Omega_{\alpha}}\Delta f<-\alpha.$$
Then for any $n\geq 1$,
$$\Omega_{\frac{\alpha}{n}}\subseteq \{x\in V:
\text{\rm{Deg}}_{n,\infty}(x)> n\}.$$ As a consequence, $(V, b,
\mu)$ has unbounded global weighted degree for any parameter $n\geq
1$.
\end{thm}
\begin{proof}
 In the proof of part (3) of Lemma \ref{comparison-lem}, we already
showed that for $x\in\Omega_{\frac{\alpha}{n}}$,
$$l=\frac{1}{\mu(x)}\sum_{y:f(y)>f(x)}b(x,y)>n.$$
We claim that for all $x\in \Omega_{\frac{\alpha}{n}}$, $k\in
\mathbb{N}$,
$$n<l\leq\text{\rm{Deg}}_{n, k}(x).$$
Assuming the claim, we see that for any $x\in
\Omega_{\frac{\alpha}{n}}$,
$$n<l\leq\text{\rm{Deg}}_{n, \infty}(x).$$
Hence
$$\Omega_{\frac{\alpha}{n}}\subseteq \{x\in V: \text{\rm{Deg}}_{n,\infty}(x)> n\}.$$

Now we complete the proof of the claim. For all $x\in
\Omega_{\frac{\alpha}{n}}$,
$$\text{\rm{Deg}}_{n, 0}(x)=\frac{1}{\mu(x)}\sum_{y}b(x,y)\geq\frac{1}{\mu(x)}\sum_{y:f(y)>f(x)}b(x,y)=l.$$
Assume that the claim is true for $k$. In other words, for all $x\in
\Omega_{\frac{\alpha}{n}}$,
$$\text{\rm{Deg}}_{n,
k}(x)\geq l>n.$$ Note that if $f(y)>f(x)$ for $x\in
\Omega_{\frac{\alpha}{n}}$, $y$ is necessarily in
$\Omega_{\frac{\alpha}{n}}$ and consequently,
$$\text{\rm{Deg}}_{n, k}(y)\geq l>n.$$ Thus we have
$$\text{\rm{Deg}}_{n, k+1}(x)=\frac{1}{\mu(x)}\sum_{y,\text{\rm{Deg}}_{n, k}(y)>n}b(x,y)\geq\frac{1}{\mu(x)}\sum_{y:f(y)>f(x)}b(x,y)=l$$
for any $x\in \Omega_{\frac{\alpha}{n}}$. The claim follows by
induction.

By Lemma \ref{degm<degn}, we see that for $m>n\geq 1$,
$$\Omega_{\frac{\alpha}{m}}\subseteq
\{x\in V: \text{\rm{Deg}}_{m,\infty}(x)> m\}\subseteq\{x\in V:
\text{\rm{Deg}}_{n,\infty}(x)> m\}.$$ The set
$\Omega_{\frac{\alpha}{m}}$ is nonempty for any $m>n$, so that the
function $\text{\rm{Deg}}_{n,\infty}(x)$ is necessarily unbounded
for any $n\geq 1$.
\end{proof}

\section{Khas'minskii's criterion}

Now we are ready to prove the following analogue of Khas'minskii's
criterion for stochastic completeness.
\begin{thm} \label{Khas-thm} Assume that the weighted degree function $\text{\rm{Deg}}(x)$ is unbounded for the graph $(V, b, \mu)$.
If there exists a nonnegative function $\gamma \in F$ on $V$ such
that
\begin{equation}\label{gamma1}
\gamma(x) \rightarrow +\infty \text{~~~~as~~~~}
\text{\rm{Deg}}(x)\rightarrow +\infty
\end{equation}
 and
 \begin{equation}\label{gamma2}
\Delta\gamma(x)+\lambda\gamma(x)\geq 0 \text{~~outside a set A of
bounded weighted degree}
\end{equation}
for some $\lambda>0$, then $V$ is stochastically complete.
\end{thm}

\begin{proof}
  We only need prove that $V$ satisfies the weak Omori-Yau
maximum principle. If not, there should exist a nonnegative function
$f$ on $V$ with $f^* <+\infty$ and some $\alpha>0$ such that
$$\sup_{\Omega_{\alpha}}\Delta f<-\alpha.$$

Let $$M=\sup\{\text{\rm{Deg}}(x):x\in A\}<+\infty.$$ By Lemma
\ref{comparison-lem}, changing $\alpha$ if necessary, we can assume
that $\text{\rm{Deg}}(x)>M$ for all $ x\in \Omega_{\alpha}$.

Let $$u=f-c\gamma,$$ where the parameter $c>0$ will be chosen later.

Since $f^* <+\infty$ and $$\gamma(x)\rightarrow +\infty
\text{~~~~as~~~~} \text{\rm{Deg}}(x)\rightarrow +\infty,$$ there
exists $N(c)>M$ such that
$$\sup_{\{x\in V: \text{\rm{Deg}}(x)< N(c)\}}u(x)=\sup_{V}u(x)<+\infty.$$
Let $0<\eta<\min(\frac{\alpha}{2},\frac{\alpha}{2\lambda})$. Since
$f$ cannot attain the value $f^* $, we can choose $\bar{x}$ such
that
$$f(\bar{x})>f^* -\frac{\eta}{2}.$$ Choose $c=c(\eta,\bar{x})>0$ small
enough to insure that $c\gamma(\bar{x})<\frac{\eta}{2}$.

For $n\in\mathbb{N}$, we can choose $x_n$ with
$\text{\rm{Deg}}(x_n)<N(c)$ such that $u(x_n)>u^* -\frac{1}{n}$. We
have
$$f(x_n )+\frac{1}{n}>f(x_n )-c\gamma(x_n )+\frac{1}{n}> f(\bar{x})-c\gamma(\bar{x})>f^* -\eta,$$ and
$$c\gamma(x_n )< f(x_n )- f^* +\eta+\frac{1}{n}<\eta+\frac{1}{n}.$$
So for every index $n>\frac{2}{\eta}$, $$f(x_n)>f^*
-\frac{3}{2}\eta>f^* -\alpha,$$
$$c\lambda\gamma(x_n )<\frac{3}{2}\lambda\eta<\frac{3}{4}\alpha.$$
In particular, for every index $n>\frac{2}{\eta}$, $x_n\in
\Omega_{\alpha}$. It follows that for all $n>\frac{2}{\eta}$,
$$\Delta\gamma(x_n)+\lambda\gamma(x_n)\geq 0,$$ and$$\Delta
f(x_n)<-\alpha.$$ Then
\begin{equation}\label{star1}
\begin{split}
\Delta(f-c\gamma)(x_n )
 &= \Delta f(x_n) -c\Delta\gamma(x_n )\\
 &<-\alpha + c\lambda \gamma(x_n )<-\alpha/4.
 \end{split}
 \end{equation}

On the other hand, we have
\begin{equation*}
\begin{split}
\Delta(f-c\gamma)(x_n )
 &= \Delta u(x_n )\\
 &=\frac{1}{\mu(x_n)}\sum_{y}b(x_n,y)(u(x_n)-u(y))\\
 &>-\frac{\text{\rm{Deg}}(x_n)}{n}>-\frac{N(c)}{n}.
 \end{split}
 \end{equation*}

Choosing sufficiently large $n$, we obtain a contradiction with (\ref{star1}).
\end{proof}

\begin{rem}\rm{
  Note that unlike in the case of manifolds we do not require that the
  exceptional set $A$ be compact.}
\end{rem}

A convenient version of Khas'minskii's criterion on manifolds is
given in \cite{PRSSurvey} . We give the discrete analogue here.

\begin{thm}\label{variant-Khas-thm}
If there exists a nonnegative function $\sigma\in F$ on $V$ with
$$\sigma(x)\rightarrow +\infty \text{~~~~as~~~~} \text{\rm{Deg}}(x)\rightarrow
+\infty$$ satisfying:
$$\Delta\sigma(x)+f(\sigma(x))\geq 0 \text{~~outside a set A of bounded weighted degree}$$
for some positive, increasing function $f\in C^1 ([0,+\infty))$ with
$$\int_{0}^{+\infty}\frac{dr}{f(r)}=+\infty,$$ then $V$ is
stochastically complete. 
\end{thm}
\begin{proof}
Let $$\phi(r)=\exp(\int_{0}^{r}\frac{ds}{f(s)+s}),$$ we have
$\phi(r)\rightarrow
  +\infty$ as $r\rightarrow +\infty$.

  The function $\phi(r)$ is increasing and concave since:
\begin{enumerate}
  \item $\phi'(r)=\frac{\phi(r)}{f(r)+r}>0$;
  \item $\phi''(r)=-\frac{\phi(r)f'(r)}{(f(r)+r)^2}\leq 0$.
\end{enumerate}

Therefore for $r,s \geq 0$ we have
$$\phi(r)-\phi(s)\geq\phi'(r)(r-s).$$
Thus
\begin{equation}
  \begin{split}
    \Delta\phi(\sigma(x))&=\frac{1}{\mu(x)}\sum_{y\in
    V}b(x,y)(\phi(\sigma(x))-\phi(\sigma(y)))\\
&\geq\phi'(\sigma(x))\frac{1}{\mu(x)}\sum_{y\in
V}b(x,y)(\sigma(x)-\sigma(y))\\&=\phi'(\sigma(x))\Delta\sigma(x),
  \end{split}
\end{equation}
which also shows that $\phi(\sigma(x))\in F$.
 Now,
consider $\gamma(x)=\phi(\sigma(x))$, then
$$\gamma(x) \rightarrow +\infty\text{~~~~as~~~~}\text{\rm{Deg}}(x)\rightarrow
+\infty.$$

On the complement of $A$ we have
\begin{equation}
\begin{split}
\Delta\gamma(x)+\gamma(x)&= \Delta\phi(\sigma(x))+\phi(\sigma(x)) \\
 &\geq \phi'(\sigma(x))\Delta\sigma(x)+\phi(\sigma(x))\\
 &=\phi'(\sigma(x))(\Delta\sigma(x)+\frac{\phi(\sigma(x))}{\phi'(\sigma(x))})\\
 &=\phi'(\sigma(x))(\Delta\sigma(x)+f(\sigma(x))+\sigma(x))\\
 &\geq\phi'(\sigma(x))(\Delta\sigma(x)+f(\sigma(x))\geq 0
 \end{split}
 \end{equation}
Theorem \ref{Khas-thm} applied to $\gamma(x)$ with $\lambda=1$
implies stochastic completeness.
\end{proof}

In the previous proof, we have made use of the following elementary
fact.
\begin{lem}
  Let $f\in C^1 ([0,+\infty))$ be a positive, increasing function.
  Assume further that
  $$\int_{0}^{+\infty}\frac{dr}{f(r)}=+\infty.$$
  Then
  $$\int_{0}^{+\infty}\frac{dr}{f(r)+r}=+\infty.$$
\end{lem}
For the sake of completeness, we give a proof here.
\begin{proof}
Note that the integral is only improper at $+\infty$ since $f$ is
positive and increasing on $[0,+\infty)$. Assume that the assertion
is not true for a while, we see that
$$\int_{0}^{+\infty}\frac{dr}{f(r)+r}<+\infty.$$
However, for all $x>0$, we have
$$0<\frac{x}{2}\cdot\frac{1}{f(x)+x}\leq\int_{\frac{x}{2}}^{x}\frac{dr}{f(r)+r}\leq\int_{\frac{x}{2}}^{+\infty}\frac{dr}{f(r)+r}.$$
The third integral necessarily goes to $0$ as $x$ approaches
$+\infty$. Thus there exists $r_0>0$ such that for any $r>r_0$,
$$\frac{r}{f(r)+r}\leq\frac{1}{2}.$$ It follows that $f(r)\geq r$
for all $r>r_0$. But then
$$\int_{r_0}^{+\infty}\frac{dr}{f(r)+r}\geq\int_{r_0}^{+\infty}\frac{dr}{2f(r)}=+\infty.$$
A contradiction.
\end{proof}

\section{Stability results}
In this section we show that after certain surgeries, a
stochastically incomplete graph will remain stochastically
incomplete. The weak Omori-Yau maximum principle allows us to pass
from the stability of existence of certain functions to the
stability of stochastic incompleteness. Roughly speaking, part (3)
of Lemma \ref{comparison-lem} implies that a perturbation of bounded
weighted degree does not affect the stochastically incompleteness.
This intuition is made explicit by the following theorems.
\begin{thm}\label{stability}
   Let $(V, b, \mu)$ be a graph and $W\subseteq V$.
$(W,b|_{W\times W},\mu|_{W})$ forms a subgraph. Assume that $W$ is
stochastically incomplete. If one of the following two conditions
holds, $V$ is also stochastically incomplete.
\begin{enumerate}
  \item For some $n\geq 1$, $\sup \{\text{\rm{Deg}}_{W}(x): x\in W, \exists y\in V\backslash W, b(x,y)>0\}<n$;
 \item There exists $n\geq 1$, such that $\forall x\in W,$ $$ \frac{1}{\mu(x)}\sum_{y\in V\backslash W}b(x,y)<n.$$
\end{enumerate}

\end{thm}
\begin{proof}
  (1)Since $W$ is stochastically incomplete there exists a nonnegative function $f$ on $W$ and $\alpha>0$ such that
  $$\sup_{\Omega^{W}_{\alpha}}\Delta^{W}
f<-\alpha.$$
Here
$$\Omega^{W}_{\alpha}=\{x\in W: f(x)>f^*-\alpha\},$$ and $$\Delta^{W}
f(x)=\frac{1}{\mu(x)}\sum_{y\in W}b(x,y)(f(x)-f(y))$$ for $x\in W$.

Define a function $u$ on $V$ by

\makeatletter
\let\@@@alph\@alph
\def\@alph#1{\ifcase#1\or \or $'$\or $''$\fi}\makeatother
\begin{subnumcases}
{u(x)=} (f(x)+\frac{\alpha}{n}-f^* )_+ , &$x\in W$, \label{eq:a1}\\
0, &$x\in V\backslash W$.\label{eq:a2}
\end{subnumcases}
\makeatletter\let\@alph\@@@alph\makeatother

We see that $u^* = \frac{\alpha}{n}$ and
$$\Omega^{V}_{\frac{\alpha}{n}}=\{x\in V: u(x)>0\}=\{x\in W: f(x)>f^*
-\frac{\alpha}{n}\}\subseteq \{x\in W: \text{\rm{Deg}}_{W}(x)>n\}$$
by (3) of Lemma \ref{comparison-lem}.

Thus for $x\in\Omega^{V}_{\frac{\alpha}{n}}$, $y\in V\backslash W$,
we have $b(x,y)=0$. Hence for every
$x\in\Omega^{V}_{\frac{\alpha}{n}}$
\begin{equation}
\begin{split}
\Delta^{V}u(x)&=\frac{1}{\mu(x)}\sum_{y\in
V}b(x,y)(u(x)-u(y))\\&=\frac{1}{\mu(x)}\sum_{y\in
W}b(x,y)(u(x)-u(y))\\&\leq\frac{1}{\mu(x)}\sum_{y\in
W}b(x,y)(f(x)-f(y))\\&=\Delta^{W}f(x)<-\alpha.
 \end{split}
 \end{equation}

The stochastic incompleteness of $V$ then follows from Theorem
\ref{equiv-thm}.

(2)As in (1), there's a nonnegative function $f$ on $W$ and
$\alpha>0$ such that
  $$\sup_{\Omega^{W}_{\alpha}}\Delta^{W}
f<-\alpha$$ since $W$ is stochastically incomplete by assumption.

Define a function $u$ on $V$ by

\makeatletter
\let\@@@alph\@alph
\def\@alph#1{\ifcase#1\or \or $'$\or $''$\fi}\makeatother
\begin{subnumcases}
{u(x)=} (f(x)+\frac{\alpha}{2n}-f^* )_+ , &$x\in W$, \label{eq:a1}\\
0, &$x\in V\backslash W$.\label{eq:a2}
\end{subnumcases}
\makeatletter\let\@alph\@@@alph\makeatother

We see that $u^* = \frac{\alpha}{2n}$ and
$$\Omega^{V}_{\frac{\alpha}{2n}}=\{x\in V: u(x)>0\}=\{x\in W: f(x)>f^*
-\frac{\alpha}{2n}\}.$$

So for $x\in\Omega^{V}_{\frac{\alpha}{2n}}$
\begin{equation}
\begin{split}
\Delta^{V}u(x)&=\frac{1}{\mu(x)}\sum_{y\in
V}b(x,y)(u(x)-u(y))\\&=\frac{1}{\mu(x)}\sum_{y\in
W}b(x,y)(u(x)-u(y))+\frac{1}{\mu(x)}\sum_{y\in V\backslash
W}b(x,y)(u(x)-u(y))\\&\leq\frac{1}{\mu(x)}\sum_{y\in
W}b(x,y)(f(x)-f(y))+\frac{1}{\mu(x)}\sum_{y\in V\backslash
W}b(x,y)\frac{\alpha}{2n}\\&\leq\Delta^{W}f(x)+\frac{\alpha}{2}<-\frac{\alpha}{2}.
 \end{split}
 \end{equation}
 The stochastic incompleteness of $V$ then follows from Theorem
\ref{equiv-thm}.
\end{proof}
\begin{rem}\rm{
 Part (2) of Theorem \ref{stability} was first proved by Keller and Lenz \cite{Keller-Lenz}. Our proof here is
 more elementary.}
\end{rem}

In Theorem \ref{stability} we derive stochastic incompleteness of
graphs from that of subgraphs. The weak maximum principle allows
also to obtain implications in the opposite direction, as in the
next statement.
\begin{thm}\label{stability-sub}
Let $(V, b, \mu)$ be a stochastically incomplete graph and $n\geq
1$. The subgraph $$W=\{x\in V: \text{\rm{Deg}}(x)>n\}$$ with weights
$(b|_{W\times W}, \mu|_{W})$ is stochastically incomplete as
well.\end{thm}

\begin{proof}
  There exists a nonnegative function $f$ on $V$ and
$\alpha>0$ such that $$\sup_{\Omega^{V}_{\alpha}}\Delta^{V}
f<-\alpha.$$ We will show that $f|_{W}$ is a function violating the
weak maximum principle.

From Lemma \ref{comparison-lem}, we see that $$\sup_W f=\sup_V f,$$
and $$\Omega^{W}_{\frac{\alpha}{n}}=\Omega^{V}_{\frac{\alpha}{n}}.$$

We claim that for any $x\in\Omega^{W}_{\frac{\alpha}{n}},$
$$\Delta^{W} f(x)\leq \Delta^{V} f(x)<-\alpha.$$

In fact, for $x\in \Omega^{W}_{\frac{\alpha}{n}}, y\in V\backslash
W$, we claim that $f(y)\leq f(x)$. If not $$f(y)>f(x)>f^*
-\frac{\alpha}{n},$$ so that
$y\in\Omega^{W}_{\frac{\alpha}{n}}\subseteq W$, a contradiction.

Then for any $x\in\Omega^{W}_{\frac{\alpha}{n}}$, we obtain
\begin{equation*}
\begin{split}
-\alpha>\Delta^V f(x)&=\frac{1}{\mu(x)}\sum_{y\in
V}b(x,y)(f(x)-f(y))\\&=\frac{1}{\mu(x)}\sum_{y\in
W}b(x,y)(f(x)-f(y))\\&+\frac{1}{\mu(x)}\sum_{y\in V\backslash
W}b(x,y)(f(x)-f(y))\\&\geq\frac{1}{\mu(x)}\sum_{y\in
W}b(x,y)(f(x)-f(y))=\Delta^{W} f(x).
 \end{split}
 \end{equation*}
The stochastic incompleteness of $V$ then follows from Theorem
\ref{equiv-thm}.
\end{proof}

\section{Applications to the physical Laplacian}
In this section, we apply the weak Omori-Yau maximum principle and
Khas'minskii's criterion to the physical Laplacian on an
(un-weighted) graph. We assume that $(V,E)$ is a locally finite,
connected infinite graph without loops and multi-edges where $V$ is
the set of vertices and $E$ is the set of edges. This corresponds to
the special case that $b(x,y)\in\{0,1\}, \mu(x)\equiv 1$.

As before, we use $V$ to denote the graph if no confusion arises. We
write $y\sim x$ if there's an edge connecting $x$ and $y$. In this
case, we call the vertices $x$ and $y$ neighbors. Then the weighted
degree function
$$\text{\rm{Deg}}(x)=\sum_{y\in V}b(x,y)=\#\{y\in V: y\sim x\},$$ is exactly the
number of neighbors of $x$ in $V$, i.e. $\deg (x)$.

Let $d$ be the graph metric on $V$, that is, for any two vertices
$x, y\in V$, $d(x,y)$ is the smallest number of edges in a chain of
edges connecting $x$ and $y$. We fix a point $x^*\in V$ as a root of
the graph and define
$$r(x)=d(x, x^* ).$$ A key feature of the graph metric is that if $x\sim
y$, then
$$|r(x)-r(y)|\leq 1.$$ We use further the notations $$S_{R}=\{y\in
V:r(y)=R\},$$$$B_{R}=\cup_{n=0}^{R}S_{n}=\{y\in V:r(y)\leq R\},$$
$$m_{\pm} (x)=\#\{y: y\sim x, r(y)=r(x)\pm 1\},$$ $$K_{\pm}
(r)=\max_{x\in S_r}m_{\pm} (x),$$ and $$k_{\pm} (r)=\min_{x\in
S_r}m_{\pm} (x).$$

The \emph{formal} Laplacian in this case is
\begin{equation}
  \Delta f(x)=\sum_{y,y\sim x}(f(x)-f(y)).
\end{equation}
Here $f$ can now be an arbitrary function on $V$ because of the
local finiteness. For example,
\begin{equation}
  \Delta r(x)=m_{-}(x)-m_{+}(x).
\end{equation}

The machinery of weak Omori-Yau maximum principle and Khas'minskii's
criterion can be applied in two ways.
\begin{enumerate}
  \item Choose a series $\sum_{n=0}^{\infty}a_n$ with nonnegative terms,
 and define the function $$f(x)=\sum_{n=0}^{r(x)}a_n$$ which then can be used in the weak Omori-Yau maximum principle and Khas'minskii's
criterion. Choosing the series appropriately we obtain sufficient
conditions for stochastic completeness and incompleteness.
 \item Alternatively, one can determine ``natural" values of $a_n$ by solving certain difference
equations or inequalities.
\end{enumerate}

Before going into details we would like to point out that for a
locally finite graph of unbounded degree, $\deg(x)\rightarrow
+\infty$ implies $r(x)\rightarrow +\infty$. Thus Theorem
\ref{Khas-thm} can be restated in a weaker form:

\begin{thm}\label{weak-Khas} Assume the degree function $\deg(x)$ is unbounded for the locally finite graph $(V,E)$. If there exists a nonnegative function $\gamma$ on $V$
with $$\gamma(x) \rightarrow +\infty \text{~~~~as~~~~}
r(x)\rightarrow +\infty$$ satisfying
$$\Delta\gamma(x)+\lambda\gamma(x)\geq 0 \text{~~outside a finite set A}$$
for some $\lambda>0$, then $V$ is stochastically complete.
\end{thm}

\begin{rem}\rm{
 Wojciechowski and Keller \cite{Keller-Woj} also obtained independently
  this form of Khas'minskii's criterion using a different method.}
\end{rem}


\subsection{Criterions for stochastic completeness}

In what follows, $\sum_{n=0}^{\infty}a_n$ is a series with
nonnegative terms.

\begin{thm}\label{series-thm}
If $\sum_{n=0}^{\infty}a_n =+\infty$ and for some $\lambda>0$, the
following inequality $$m_+ (x) a_{r(x)+1} -m_- (x)a_{r(x)}\leq
\lambda \sum_{n=0}^{r(x)}a_n$$ holds outside a finite set, then $V$
is stochastically complete.
\end{thm}
\begin{proof}
  Let $\gamma(x)=\sum_0^{r(x)}a_n $, then $$\Delta \gamma(x)+\lambda\gamma(x) =m_- (x)a_{r(x)}-m_+ (x)
  a_{r(x)+1}+\lambda \sum_0^{r(x)}a_n \geq 0$$
outside a finite set and $\gamma(x)\rightarrow +\infty$ as
$r(x)\rightarrow+\infty$. By Theorem \ref{weak-Khas}, $V$ is
stochastically complete.
\end{proof}

Theorem \ref{series-thm} already gives some nontrivial results
through some obvious choices of $a_n$:
\begin{enumerate}
\item One natural choice is $a_n \equiv 1$. Then a sufficient condition
for stochastic completeness is $$m_+ (x) -m_- (x)\leq \lambda r(x)$$
outside a finite set for some $\lambda>0$. This improves the
curvature type criterion of Weber \cite{Weber} where the sufficient
condition is $$m_+ (x) -m_-
(x)\leq C$$ for some constant $C>0$. 

\item Take $a_0=0,a_n =\frac{1}{n}$ for $n>1$. We have that $$m_+ (x)- (1+\frac{1}{r(x)})m_- (x)\leq \lambda \log
r(x)$$ is sufficient for stochastic completeness.
\end{enumerate}

One can improve these results by choosing divergent series with
smaller terms. We do this with a view toward using Theorem
\ref{variant-Khas-thm}.

\begin{thm}\label{curvature-thm}
  If for some positive, increasing function $f\in C^1 ([0,+\infty))$ with
$\int_{0}^{+\infty}\frac{dr}{f(r)}=+\infty$, $$m_+ (x) -m_- (x)\leq
f(r(x))$$ outside a finite set,then $V$ is stochastically complete.
\end{thm}
\begin{proof}
  We only need to take $\sigma(x)=r(x)$ in Theorem \ref{variant-Khas-thm}.
\end{proof}

\begin{rem}
\rm{The quantity $\Delta r(x)=m_- (x) -m_+ (x)$ can be viewed as an analogue of the mean curvature of a geodesic sphere on a Riemannian manifold.}
\end{rem}

The following result was first obtained by Wojciechowski
\cite{WojSurvey}. We give a shorter proof, based on Theorem
\ref{weak-Khas}.

\begin{thm}
  If $\sum_{r=0}^{\infty}\frac{1}{K_+ (r)}=+\infty$, then $V$ is
  stochastically complete.
\end{thm}
\begin{proof}
  Let $$\gamma(x)=\sum_{r=0}^{r(x)-1}\frac{1}{K_+ (r)}$$ for $r(x)>0$,
  and $\gamma(x^* )=0$. We then have that
  $$\gamma(x)\rightarrow+\infty\text{~~~~as~~~~}r(x)\rightarrow+\infty,$$
  and outside a finite set
  $$\Delta\gamma(x)+\gamma(x)=m_- (x)\frac{1}{K_+ (r(x)-1)}-m_+ (x)\frac{1}{K_+ (r(x))}+\gamma(x)\geq \gamma(x)-1 \geq 0.$$
   The assertion follows from Theorem \ref{weak-Khas}.
\end{proof}

\subsection{Criteria for stochastic incompleteness}%

Similarly, using test series to define functions that violate the
weak maximum principle, we obtain a curvature type criterion for
stochastic incompleteness:

\begin{thm} If $\sum_{l=0}^{\infty}a_l <+\infty, a_l \geq 0$ and
for some $n\in \mathbb{N}, c>0$, the inequality $$m_+ (x) a_{r(x)+1}
-m_- (x)a_{r(x)}>c$$ holds for $r(x)>n$, then $V$ is stochastically
incomplete.
\end{thm}

\begin{proof}
Let $$f(x)=\sum_{l=0}^{r(x)}a_l.$$ Then $$f^*
=\sum_{r=0}^{\infty}a_r <+\infty.$$

Let $\alpha=\sum_{l=n+1}^{\infty}a_l.$ Then $f(x)>f^* -\alpha$
implies that $r(x)>n$. So in this case,
$$-\Delta f(x)=m_+ (x) a_{r(x)+1} -m_-
(x)a_{r(x)}>c.$$ By Theorem \ref{equiv-thm}, $V$ is stochastically
incomplete.
\end{proof}

Theorem \ref{equiv-thm} can also be used to derive the following
result about stochastic incompleteness obtained by Wojciechowski
\cite{WojIndiana}.

\begin{thm}\label{ratio-curvature-thm}
  If
  \begin{equation}
    \label{m_pm}
\sum_{r=1}^{\infty}\max_{x\in
S_{r}}\frac{m_{-}(x)}{m_{+}(x)}<+\infty,
  \end{equation}
then $V$
  is stochastically incomplete.
\end{thm}
\begin{proof}
Denote $\max_{x\in S_{r}}\frac{m_{-}(x)}{m_{+}(x)}$ by $\eta(r)$.
Let $$f(x)=\sum_{r=1}^{r(x)-1}\eta(r)$$ for $r(x)\geq 2$, and
$f(x)=0$ elsewhere. Then $$f^* = \sup
f(x)=\sum_{r=1}^{\infty}\eta(r)<+\infty.$$ Choose $r_0>2$
sufficiently large so that
$$0<\alpha=\sum_{r=r_0 -1}^{\infty}\eta(r)<\frac{1}{2}.$$ Then
$$\Omega_{\alpha}=\{x\in V: f(x)>\sum_{r=1}^{r_0-2}\eta(r)\}=B_{r_0 -1}^c.$$ But for $x\in
B_{r_0 -1}^c,$
$$\eta(r(x)-1)<\frac{1}{2}.$$ Hence
\begin{equation*}
\begin{split}
\Delta
f(x)&=m_{-}(x)\eta(r(x)-1)-m_{+}(x)\eta(r(x))\\
&\leq \frac{1}{2}m_{-}(x)- m_{-}(x)\\&\leq -\frac{1}{2}m_{-}(x)\leq
-\frac{1}{2}
 \end{split}
 \end{equation*}
on $\Omega_{\alpha}$. 

By Theorem \ref{equiv-thm}, $V$ is stochastically incomplete.
\end{proof}

\begin{rem}\rm{Theorem \ref{ratio-curvature-thm} first appeared in a slightly weaker form as Theorem 3.4 in \cite{WojIndiana}.
There stochastic incompleteness is established under the condition
$$\sum_{r=1}^{\infty}\frac{K_{-}(r)}{k_{+}(r)}<+\infty$$ instead of
(\ref{m_pm}).}
\end{rem}

\subsection{The symmetric case}

\begin{defi}
A graph $V$ is called weakly symmetric if it satisfies $$m_+ (x)=g_+
(r(x)),m_- (x)=g_- (r(x))$$ with functions $g_+ (r),g_- (r):
\mathbb{N}\rightarrow \mathbb{N}$.
\end{defi}

For graphs that are weakly symmetric, Wojciechowski \cite{WojSurvey}
proved the following criteria. Here we present a proof based on the
weak maximum principle.

\begin{thm}\label{sym-thm}
 A weakly symmetric graph $V$ is stochastically complete if and only if $$\sum_{r=0}^{\infty}\frac{V(r)}{g_+
  (r)S(r)}=+\infty$$ where $S(r)=\#S_{r}$ and
  $V(r)=\#B_{r}$.
\end{thm}
\begin{proof}
Since $$m_+ (x)=g_+ (r(x)),m_- (x)=g_-
  (r(x)),$$ we see that
  $$g_- (r)S(r)=g_+ (r-1)S(r-1).$$

Let $$\gamma(x)=\sum_{r=0}^{r(x)-1}\frac{V(r)}{g_+
  (r)S(r)}$$ for $r(x)>0$,
  and $\gamma(x^* )=0$. We have
  \begin{equation*}
\begin{split}
\Delta\gamma(x)&=g_- (r(x))\frac{V(r(x)-1)}{g_+
  (r(x)-1)S(r(x)-1)} -g_+ (r(x))\frac{V(r(x))}{g_+
  (r(x))S(r(x))}\\&=\frac{V(r(x)-1)}{S(r(x))}-\frac{V(r(x))}{S(r(x))}=-1
 \end{split}
 \end{equation*}
for $r(x)\geq 1$.

If $\gamma(x)\rightarrow+\infty$ as $r(x)\rightarrow+\infty$, then
  $$\Delta\gamma(x)+\gamma(x)=\gamma(x)-1\geq 0$$ outside a finite
  set. The stochastic
incompleteness then follows from Theorem \ref{weak-Khas}.

For the other implication suppose that $\gamma^* = \sup
  \gamma(x)<+\infty$. Letting $\alpha=\gamma^* $, we see that on
$\Omega_{\alpha}=B_0^c$, $$\Delta \gamma(x)=-1.$$ The stochastic
incompleteness then follows from Theorem \ref{equiv-thm}.
\end{proof}

\begin{rem}\rm{
As pointed out by Wojciechowski \cite{WojSurvey}, it is interesting
to notice that for a weakly symmetric graph, the edges between
points on the same sphere play no role in stochastic completeness.
See also \cite{Keller-Woj} for further studies of weakly symmetric
graphs.}
\end{rem}

\section{Further remarks}


(1) A rich source for ideas behind the study of stochastic
completeness of the physical Laplacian is the literature about the
Riemannian manifold case. However, due to the fact that the
Dirichlet form on a graph is nonlocal, there are some essential
differences in our case. For example, as shown by Wojciechowski
\cite{WojSurvey}, there exist stochastically incomplete graphs with
polynomial volume growth which never happens in the manifold case.
His examples of stochastically incomplete graphs that satisfy
$$\mu(B_r)\leq Cr^{3+\varepsilon}, C,\varepsilon>0,$$ are presented
in the next remark.

It is then interesting to ask what is the smallest possible volume
growth for stochastically incomplete graphs. It is natural to
conjecture that for the physical Laplacian on graphs, the condition
$$\mu(B_r)\leq Cr^{3}, C>0,$$
implies stochastic completeness. This is proven in a forthcoming
paper of Grigor'yan, Huang and Masamune \cite{Gri-H-Mas}. Note that
for geodesically complete Riemannian manifolds, the almost sharp
condition (\cite{Davies2}, \cite{Grigor2}, \cite{Hsu},
\cite{Karp-Li})
$$\mu(B_r)\leq\exp Cr^2, C>0,$$ implies
stochastic completeness.

On the other hand, there exist stochastically complete graphs with
arbitrarily large volume growth. For example, take a set of vertices
$\{0,1,2,\cdots,n \cdots \}$ with edges $n\sim n+1$. For each vertex
$n$, we associate a distinct finite set $V_n$ and add extra edges
between $n$ and points in $V_n$. The resulting graph $V$ is then a
tree whose volume growth can be chosen to be arbitrarily large. It
is of bounded global weighted degree with parameter $1$ and hence is
stochastically complete by Theorem \ref{global-degree}. The
stochastic completeness of $V$ can be shown via Theorem
\ref{stability-sub} as well.

(2) Let $S(r)$ be given with $S(0)=1$. By connecting every vertex in
$S_{r}$ to every vertex in $S_{r+1}$ we get a spherically symmetric
graph $G_{S}$. Then by Theorem \ref{sym-thm}, $G_{S}$ is
stochastically incomplete if and only if
$$\sum_{r=0}^{\infty}\frac{\sum_{i=0}^{r}S(i)}{S(r+1)S(r)} < +\infty$$
since $m_{\pm} (x)= S(r(x)\pm1)$. Taking $S(r)=
[(r+1)^{2+\varepsilon}], \varepsilon>0$ where $[c]$ is the integer
part of $c$, we see that $G_{S}$ is stochastically incomplete
whereas
$$\mu(B_r)\leq C r^{3+\varepsilon}$$
for some $C>0$.

This construction of Wojciechowski \cite{WojSurvey}, at the same
time gives a counterexample to the converse of Theorem
\ref{ratio-curvature-thm}. The graph $G_{S}$ with $S(r)= (r+1)^{3}$
satisfies
$$\sum_{r=1}^{\infty}\max_{x\in
S_{r}}\frac{m_{-}(x)}{m_{+}(x)}=\sum_{r=1}^{\infty}\frac{(r-1)^3}{(r+1)^3}=+\infty,$$
but is stochastically incomplete.

(3) We conjecture that if:
$$m_{+}(x) -m_{-}(x) \geq f(r(x)),$$ where $f(r)>0$ and $\sum_{r=0}^{\infty}\frac{1}{f(r)}<+\infty$, then $V$ is
stochastically incomplete. This should be a useful complement to
Theorem \ref{curvature-thm}.

(4) We conjecture that the converse of Theorem \ref{Khas-thm} should
be true. Namely, if a graph $(V, b, \mu)$ is stochastically
complete, then there should exist a function $\gamma(x)\in F$ on $V$
satisfying the conditions (\ref{gamma1}), (\ref{gamma2}).

(5) For a subset $A$ of $V$, we define its (outer) boundary to be
$$\partial A =\{x: x\in A^c ,\text{~~~~and~~~~} \exists y\in A, \text{~~~~s.t.~~~~} x\sim y
\}$$ and its closure to be $\bar{A}= A\cup\partial A$. Wojciechowski
and Keller \cite{Keller-Woj} proposed the following conjecture.

\begin{conj}
If for some fixed point $x^* \in V$ as root,
\begin{equation}\label{KL-conj}
\sum_{r=0}^{\infty}\frac{\#B_r}{\#\partial B_r} =+\infty,
\end{equation}
 then
$(V, E)$ is stochastically complete.
\end{conj}

This is an analogue of a conjecture for the stochastic completeness
of manifold proposed by Grigor'yan in \cite{Grigor}. However,
recently B\"{a}r and Bessa \cite{Bar-Bessa} constructed a
counterexample to Grigor'yan's conjecture. Their idea can also be
applied to the physical Laplacian as follows.

Take a stochastically complete tree $T$ with root $x_1$, for
example, a binary tree. Then $T$ has exponential volume growth with
respect to graph distance. Choose a stochastically incomplete graph
with only polynomial volume growth, for example, the graph $G_S$ in
the previous remark with the root denoted by $x_2$. Now we make a
single extra edge between $x_1$ and $x_2$ resulting in a new graph
$V$. Since the gluing happens at only one point at $G_S$, the graph
$V$ is stochastically incomplete by Theorem \ref{stability}.
However, for any fixed point $x^* \in V$ as a root, the quantities
$\#B_r$ and $\#\partial B_r$ are always of the order $2^n$. So we
know that $V$ satisfies (\ref{KL-conj}) while it is stochastically
incomplete. This example is simpler than the example of
\cite{Bar-Bessa} in the manifold case, thanks to special features of
the discrete setting.

\section*{Acknowledgement}
The author is grateful to his supervisor, Prof. Grigor'yan, who
introduced this topic to him and made several valuable suggestions.
The author also would like to thank R. Wojciechowski, M. Keller, and
D. Lenz for inspiring discussions and generously sharing of
knowledge. Part of this work was done when the author was visiting
the University of Jena. He would like to thank the Department of
Mathematics and Computer Science there for its hospitality.

\end{document}